\input amstex
\documentstyle{amsppt}
\magnification1200

\def\w{\omega}
\def\IN{\Bbb N}
\def\IZ{\Bbb Z}
\def\e{\varepsilon}
\def\ti{\times}
\def\K{\Cal K}
\NoBlackBoxes

\vsize=18cm

\topmatter
\title
Topologies on groups determined by sequences:\\
Answers to several questions of I.Protasov and E.Zelenyuk
\endtitle
\rightheadtext{Topologies on groups determined by sequences}
\author
Taras Banakh
\endauthor
\address
Department of Mechanics and Mathematics, Lviv Ivan Franko National University,
Universytetska 1, Lviv, 79000, Ukraine
\endaddress
\email
tbanakh\@franko.lviv.ua
\endemail
\thanks Research supported in part by grant INTAS-96-0753.
\endthanks
\abstract
Answering questions of Protasov and Zelenyuk we prove the
following results:
\item{1.} For every
increasing function $f:\IN\to\IN$ with $\lim_{n\to\infty} f(n+1)-f(n)=\infty$
and every metrizable totally bounded group topology $\tau$ on $\IZ$
there exists a convergent to zero sequence $(a_n)_{n\in\w}$ in $(\IZ,\tau)$
such that $\lim_{n\to\infty}\frac{a_n}{f(n)}=1$.
\item{2.} For every real $r>1$ there exists a sequence
$(a_n)_{n\in\w}\subset\IZ$ such that
$\lim_{n\to\infty}\frac{a_{n+1}}{a_n}=r$ but there is no ring topology
$\tau$ on $\IZ$ such that $(a_n)_{n\in\w}$ converges to zero in
$(\IZ,\tau)$.
\item{3.} There exists a countable topological Abelian group $G$ determined
by a $T$-sequence and containing a closed subgroup $H$ which is not determined
by a $T$-sequence but is homeomorphic to $G$.
\item{4.} There exist two group topologies $\tau_1$, $\tau_2$ determined by
$T$-sequences on $\IZ$ such that the topology $\tau_1\lor\tau_2$ is not
complete and thus is not determined by a $T$-sequence.
\item{5.} There exists a countable Abelian group admitting a group topology $\tau$
determined by a $T$-sequence and a metrizable group topology  $\tau'$ such that
the topology $\tau\lor\tau'$ is not discrete but contains no non-trivial
convergent sequence.

\endabstract
\endtopmatter

\document
In this note we give answers to several problems posed by I.Protasov and
E.Zelenyuk in \cite{PZ$_1$} and \cite{PZ$_2$}.
Following \cite{PZ$_2$} we define a sequence $(a_n)_{n\in\w}$ of
elements of
a group $G$ to be a {\it $T$-sequence} if $(a_n)_{n\in\w}$ converges to zero in
some non-discrete Hausdorff group topology on $G$. Given a $T$-sequence
$(a_n)_{n\in\w}$ in $G$ we denote by $(G|(a_n))$ the
group $G$ endowed with the strongest topology in which the sequence $(a_n)$
converges to zero. We say that a topological group $G$ is {\it determined
by a $T$-sequence} if $G=(G|(a_n))$ for some $T$-sequence
$(a_n)_{n\in\w}$ in $G$.

\heading
1. There is no restriction on the growth of $T$-sequences in $\IZ$
\endheading

All $T$-sequences of integers constructed in \cite{PZ$_1$} have
exponential growth. This led I.Protasov and E.Zelenyuk to the following
question (see \cite{PZ$_2$} and \cite{PZ$_1$, Question 2.2.3}): {\sl is
there a monotone $T$-sequence of integers having polynomial growth}?
First our result answers this question affirmatively. We recall that a
group topology $\tau$ on a group $G$ is called {\it totally bounded} if
for every neighborhood $U\in\tau$ of zero in $G$ there exists a finite
subset $F\subset G$ with $G=F\cdot U$.

\proclaim{Theorem 1}
\roster
\item If $(a_n)_{n=1}^\infty\subset\IZ$
is an increasing $T$-sequence, then
$\lim_{n\to\infty}(a_{n+1}-a_n)=\infty$.
\item Suppose $f:\IN\to\IN$ and $\e:\IN\to[0,\infty)$ are functions such
that $\lim_{n\to\infty}\e(n)=\infty$ and
$\lim_{n\to\infty}f(n+1)-f(n)=\infty$.
For every metrizable totally bounded
group topology $\tau$ on $\IZ$ there exists a converging to zero sequence
$(a_n)_{n\in\w}\subset(\IZ,\tau)$ such that
$\lim_{n\to\infty}\frac{a_n}{f(n)}=1$ and $|a_n-f(n)|\le\e(n)$ for every $n\in\w$.
\endroster
\endproclaim

\demo{Proof} 1. Suppose $(a_n)_{n\in\w}\subset\IZ$ is an increasing
$T$-sequence with
$\lim_{n\to\infty}a_{n+1}-a_n\ne\infty$. This means that for some
$C\in\IN$ and every $n\in\IN$ we can find $m\ge n$ with
$a_{m+1}-a_m\le C$. Let $\tau$ be a non-discrete Hausdorff group topology on $\IZ$
such that $(a_n)_{n=1}^\infty$ converges to zero in $\tau$. Pick a
$\tau$-open neighborhood $U\subset\IZ$ of zero such that $U\cap
(i+U)=\emptyset$ for every $1\le i\le C$ and find $n_0\in\IN$ such that
$a_n\in U$ for every $n\ge n_0$. By the choice of the constant $C$, there
exists $m\ge n_0$ with $a_{m+1}-a_m\le C$. Then letting $i=a_{m+1}-a_m$, we
get $a_{m+1}=a_m+i\in (i+U)\cap U=\emptyset$, a contradiction.

2. Suppose functions $f$ and $\e$ satisfy the hypotheses of the theorem.
Without loss of generality, $\e(1)=0$ and
$\e(n)\le\frac12\min\{\sqrt{f(n)},f(n+1)-f(n),f(n)-f(n-1)\}$ for $n>1$.

Let $\tau$ be any metrizable totally bounded group topology on $\IZ$ and
$\IZ=U_0\supset U_1\supset U_2\supset\dots$ be a countable base of
neighborhoods of zero in $(\IZ,\tau)$. For every $n\in\w$ let
$k(n)=\max\{i\in\w:U_i\cap[f(n)-\e(n),f(n)+\e(n)]\ne\varnothing\}$
and $a_n$ be any point in $U_{k(n)}\cap[f(n)-\e(n),f(n)+\e(n)]$ (the
number $k(n)$ is finite since the topology $\tau$ is Hausdorff).
Evidently, $|f(n)-a_n|\le\e(n)$ for every $n\in\w$ and
$0\le\lim_{n\to\infty}\big|\frac{a_n}{f(n)}-1\big|\le
\lim_{n\to\infty}\frac{\e(n)}{f(n)}\le\lim_{n\to\infty}\frac1{\sqrt{f(n)}}=0$.

It remains to verify the convergence of the constructed sequence
$(a_n)_{n\in\w}$ to zero in the topology $\tau$. This will follow as soon
as we prove that $\lim_{n\to\infty}k(n)=\infty$. Fix any number $m\in\IN$.
We have to find $n_0\in\IN$ such that $k(n)\ge m$ for every $n\ge n_0$.
Using the total boundedness of the topology $\tau$, find
$l\in\IN$ such that $\bigcup_{|i|<l}(i+U_m)=\IZ$. Since
$\lim_{n\to\infty}\e(n)=\infty$, there exists $n_0\in\IN$ such that
$\e(n)>l$ for all $n\ge n_0$. It follows that for every $n\ge n_0$ there
exists $i\in\IZ$ such that $|i|<l<\e(n)$ and $i+U_m\ni f(n)$. Consequently,
$U_m\cap[f(n)-\e(n),f(n)+\e(n)]\ne\varnothing$ and hence $k(n)\ge m$.
\qed
\enddemo

\remark{Remark 1} The requirement of the metrizability of the topology
$\tau$ in Theorem 1 is essential: according to \cite{PZ$_1$, \S5.1}, there
exists a totally bounded group topology $\tau$ on $\IZ$ such that the
space $(\IZ,\tau)$ contains no nontrivial convergent sequence.
\endremark

\remark{Remark 2} Theorem 1 gives a short proof of Theorem 2.2.6 from
\cite{PZ$_1$} which states that for every real number $r\ge 1$ there
exists a $T$-sequence $(a_n)_{n\in\w}\subset\IZ$ with
$\lim_{n\to\infty}\frac{a_{n+1}}{a_n}=r$ (apply Theorem 1 with
$f(n)=r^n+n^2$ and $\e(n)=n$).
\endremark

\heading
2. $T$-sequences in the ring $\IZ$.
\endheading

According to Theorem 2.2.3 \cite{PZ$_1$}, if $(a_n)_{n\in\w}\subset\IZ$ is
a sequence such that $\lim_{n\to\infty}\frac{a_{n+1}}{a_n}$ is a
transcendental real number, then $(a_n)_{n\in\w}$ is a $T$-sequence in the
group $\IZ$. In \cite{PZ$_2$} (see also \cite{PZ$_1$, Question 3.4.1})
I.Protasov and E.Zelenyuk asked if such a sequence $(a_n)_{n\in\w}$ needs
to be a $T$-sequence in the ring $\IZ$, i.e., if $(a_n)_{n\in\w}$
converges to zero for some Hausdorff ring topology $\tau$ on $\IZ$.
The following theorem answers  this question in negative.

\proclaim{Theorem 2} For every real number $r>1$ there exists a sequence
$(a_n)_{n\in\w}\subset\IZ$ such that
$\lim_{n\to\infty}\frac{a_{n+1}}{a_n}=r$ but $(a_n)_{n\in\w}$ is not a
$T$-sequence in the ring $\IZ$.
\endproclaim

\demo{Proof} Given a real number $r>1$ consider the sequence
$(a_n)_{n\in\w}\subset\IZ$ defined by
$$
a_n=\cases
[r^{n/2}]^2+1,&\text{ if $n=2\cdot 3^k$ for some $k\in\IN$};\\
[r^n],& \text{ otherwise},
\endcases
$$
where as usual $[x]=\max\{k\in\IZ: k\le x\}$ for a real number $x$. It can
be easily shown that $\lim_{n\to\infty}\frac{a_{n+1}}{a_n}=r$.
Nonetheless, the sequence $(a_n)_{n\in\w}$ can not converge to zero in a
ring topology $\tau$ on $\IZ$ because $a_{2\cdot 3^k}-a^2_{3^k}=1$ for
every $k\in\w$.\qed
\enddemo

\heading 3. On closed subgroups of topological groups determined by
$T$-sequences.
\endheading

In this section we give an example of a countable Abelian topological
group $G$ determined by a $T$-sequence and a closed subgroup $H$ of $G$
which is not determined by a $T$-sequence, thus answering Question 2.3.1
of \cite{PZ$_1$}. The group $G$ is a
Graev free topological Abelian group $A(S_0)$ over a convergent sequence
$S_0$ under which we understand any countable compactum $S_0$ with a unique
nonisolated point considered as the distinguished point of $S_0$.
We recall that {\it the Graev free topological Abelian group} $A(X)$ over a
pointed Tychonov space $(X,*)$ is
uniquely determined by the following three requirements:
(1) there is an embedding $X\subset A(X)$ such that the fixed point $*$ of $X$
coincides with the neutral element of the group $A(X)$,
(2) $A(X)$ coincides with the group hull of $X$ in $A(X)$, and (3) every
continuous map $f:X\to G$ into a topological Abelian group $G$ such that
$f(*)=0$ uniquely extends to a continuous group homomorphism
$\bar f:A(X)\to G$, see \cite{Gr}.

It will be more convenient to work with the following
concrete realization of a free group $A(S_0)$. In the group $\IZ^\w$
consider the sequence $(e_n)_{n\in\w}\subset\IZ^\w$ of characteristic
functions $e_n=\chi_{\{n\}}:\w\to\{0,1\}\subset\IZ$ of one-point subsets
$\{n\}\subset\w$. Clearly, the sequence $(e_n)_{n\in\w}$ converges to zero
in the product topology of $\IZ^\w$. Denote by $\IZ^\w_f$ the group hull
of the set $\{e_n:n\in\w\}$ in $\IZ^\w$.
Algebraically, $\IZ^\w_f$ is the direct sum of countably many of cyclic
groups $\IZ$ and consists of all eventually zero sequences of integers.
 It can be easily shown that a
Graev free topological Abelian group $A(S_0)$ is topologically isomorphic
to the group $(\IZ^\w_f|(e_n))$ determined by the
$T$-sequence $(e_n)_{n\in\w}$.

\proclaim{Theorem 3} The topological group $A(S_0)=(\IZ^\w_f|(e_n))$ contains a closed subgroup $H$ which is not determined by a
$T$-sequence.
\endproclaim

\demo{Proof}
 Fix any function $f:\w\to\IN$ such that the set $f^{-1}(n)$ is infinite
for every $n\in\IN$ and consider the closed subgroup
$$
H=\big\{(x_i)_{i\in\w}\in \IZ^\w_f:x_i\in f(i)\cdot\IZ\text{ \ for every \
}i\in \w\big\}
$$
in $(\IZ^\w_f|(e_n))$. We claim
that the topology of $H$ is not determined by a $T$-sequence.
Suppose the contrary: $H=(H|(b_n))$ for some sequence
$(b_n)_{n\in\w}$ convergent to zero in the topology $\tau$ of the group
$(\IZ^\w_f|(e_n))$. By the standard arguments (see Chapter 4
of \cite{PZ$_1$}) it can be shown that
the group $(\IZ^\w_f|(e_n))$ carries the strongest topology inducing the original
(product) topology on each (compact) set
$$
\IZ^\w_n=\{(x_i)_{i\in\w}\in\IZ^\w_f:\sum_{i\in\w}|x_i|\le n\},\;\; n\in\IN.
$$
This fact and the convergence of $(b_n)_{n\in\w}$ in $(\IZ^\w_f|(e_n))$ imply
$\{b_n:n\in\w\}\subset \IZ^\w_{n_0}\cap H$ for some $n_0\in\IN$.
Observe that $H=H'\oplus H''$, where $H'$ and $H''$ are the groups hulls
of the sets
$\{f(n)e_n:f(n)\le n_0\}$ and $\{f(n)e_n:f(n)>n_0\}$,
respectively. It follows from the inclusion $\{b_n\}_{n\in\w}
\subset \IZ^\w_{n_0}\cap H$ that $\{b_n\}_{n\in\w}\subset H'$.
Since $H=H'\oplus H''$ carries the strongest group
topology in which the sequence $(b_n)_{n\in\w}$ converges to zero, we
conclude that the topology of the group $H''$ is discrete. But this is not
so, because $H''$ contains an infinite compact subset
$\{0\}\cup\{(n_0+1)e_i:f(i)=n_0+1\}$. \qed
\enddemo

\remark{Remark 3} Using the Zelenyuk topological classification of
countable $k_\w$-group (\cite{Ze} or \cite{PZ$_1$, \S4.3}) it can be shown that
the subgroup $H\subset A(S_0)$ constructed in Theorem 3 is homeomorphic to
$A(S_0)$. This shows that the property of a topological group to be
determined by a $T$-sequence is not a topological invariant.
\endremark

\remark{Remark 4} The pathology described in Theorem 3 can not occur in
the group $\IZ$: for every Hausdorff group topology $\tau$ on $\IZ$ every
closed non-trivial subgroup $H$ of $(\IZ,\tau)$ has finite index in $\IZ$
and thus is open. If $(\IZ,\tau)$ is determined by a $T$-sequence, then
so does the open subgroup $H$.
\endremark

\proclaim{Question} Suppose $G$ is a finitely-generated Abelian
topological group determined by a $T$-sequence. Is every closed subgroup
of $G$ determined by a $T$-sequence?
\endproclaim

\heading
4. On supremum of group topologies determined by $T$-sequences
\endheading

In \cite{PZ$_2$} I.Protasov and
E.Zelenyuk posed the following problem: {\sl Suppose $\tau_1,\tau_2$ are two group
topologies on $\IZ$ determined by $T$-sequences.
Is the topology $\tau_1\lor \tau_2$ determined by a $T$-sequence?} We
recall that $\tau_1\lor \tau_2$ is the weakest topology $\tau$ on $\IZ$
such that the identity maps $(\IZ,\tau)\to (\IZ,\tau_i)$, $i=1,2$, are
continuous. Clearly, the group $(\IZ,\tau_1\lor\tau_2)$ may be identified
with the diagonal of the product $(\IZ\ti\IZ,\tau_1\ti\tau_2)$.

It turns out that the supremum $\tau_1\lor\tau_2$ of two topologies
determined by $T$-sequences on a countable Abelian group $G$ may be very
wild: it needs not be a $k_\w$-topology as well as may be a
$k_\w$-topology but not determined by a $T$-sequence, etc.

We remind that a topological group $G$ is called a {\it $k_\w$-group} if
$G$ admits a cover $\K$ by compact subspaces such that a subset $U\subset
G$ is open in $G$ if and only if $U\cap K$ is open in $K$ for every
$K\in\K$. According to \cite{PZ$_1$, Corollary 4.1.5} every countable group $G$
determined by a $T$-sequence is a $k_\w$-group.

Given an Abelian group $G$ and a subgroup $H$ of $G$ let $G\oplus_HG$
denote the quotient group of $G\oplus G$ by the subgroup
$\Gamma=\{(h,-h):h\in H\}\subset G\oplus G$. The following result was
suggested by I.V.Protasov.

\proclaim{Theorem 4} For every subgroup $H$ of a topological group $G$
determined by a $T$-sequence there exist group topologies $\tau_1$,
$\tau_2$ on $G\oplus_H G$ determined by $T$-sequences such that the
topological group $(G\oplus_HG,\tau_1\lor\tau_2)$ contains an open
subgroup topologically isomorphic to the group $H$.
\endproclaim

\demo{Proof} Let $(a_n)_{n\in\w}\subset G$ be a $T$-sequence determining
the topology of the group $G$. Denote by $\pi:G\oplus G\to G\oplus_HG$ the
quotient homomorphism and by $e_1,e_2:G\to G\oplus_HG$ the injective group
homomorphisms defined by $e_1(g)=\pi(g,0)=(g,0)+\Gamma$ and
$e_2(g)=\pi(0,g)=(0,g)+\Gamma$ for $g\in G$. Observe that
$e_1(h)=(h,0)+\Gamma=(0,h)+\Gamma=e_2(h)$ for any $h\in H$ which allows us
to define the injective homomorphism $e:H\to G\oplus_HG$ by
$e=e_1|H=e_2|H$. It is easy to see that $e(H)=e_1(G)\cap e_2(G)$.
Using Theorem 2.1.4 of \cite{PZ$_1$} (or the complementability of the
groups $e_i(G)$ in $G\oplus_HG$) one may show that for $i=1,2$ the
sequence $(e_i(a_n))_{n\in\w}$ is a $T$-sequence in $G\oplus_HG$
determining the non-discrete topology $\tau_i$ on $G\oplus_HG$. It
follows that the map $e_i:G\to(G\oplus_HG,\tau_i)$ is an open embedding.
Then the map $e:H\to(G\oplus_HG,\tau_1\lor\tau_2)$ is a topological
embedding and $e(H)=e_1(G)\cap e_2(G)$ is an open subgroup of
$(G\oplus_HG,\tau_1\lor\tau_2)$ isomorphic to $H$.\qed
\enddemo

Theorem 4 shows that the supremum $\tau_1\lor\tau_2$ of two group
topologies determined by $T$-sequences may be as bad as bad are subgroups
of topological groups determined by $T$-sequences.
In particular, Theorems 3 and 4 imply

\proclaim{Corollary 1} There exists a countable Abelian group $G$ and two
topologies $\tau_1$, $\tau_2$ on $G$ determined by $T$-sequences such that
$(G,\tau_1\lor\tau_2)$ is a $k_\w$-group not determined by a $T$-sequence.
\endproclaim

Theorem 4 also yields the existence
of a countable Abelian group $G$ and two
topologies $\tau_1$, $\tau_2$ on $G$ determined by $T$-sequences such that
$(G,\tau_1\lor\tau_2)$ is not complete and thus is not a $k_\w$-group.
We shall show that this fact is valid even for the group $G=\IZ$
thus answering Question 3 of \cite{PZ$_2$}
(we do not know if Corollary 1 is true for the group $G=\IZ$).

It follows from Theorems 2.2.3 and 2.2.1 of \cite{PZ$_1$} (see also
Exercise 2.2.5 in \cite{PZ$_1$}) that a sequence
$(a_n)_{n\in\w}\subset\IZ$ is a $T$-sequence in $\IZ$ provided
$\lim_{n\to\infty}\frac{a_{n+1}}{a_n}$ either is infinite or is a
transcendental real number. This implies that for any such a $T$-sequence
$(a_n)_{n\in\w}$ and any non-zero integer $c$ the sequence $(a_n+c)_{n\in\w}$ is a
$T$-sequence too. Denote by
$\tau_1$, $\tau_2$ the group topologies
on $\IZ$ determined by the $T$-sequences $(a_n)_{n\in\w}$, $(a_n+c)_{n\in\w}$.
The following theorem implies that the topology $\tau_1\lor\tau_2$ on $\IZ$
is not complete and thus is not determined by a $T$-sequence.

\proclaim{Theorem 5} Let $\tau_1$ and $\tau_2$ be the group topologies on
a countable group $G$
determined by $T$-sequences $(a_n)_{n\in\w}$ and $(b_n)_{n\in\w}$,
respectively.
If for some non-zero element $g\in G$ and every $n_0\in\w$ there exist
$n,m\ge n_0$ with $g=a_n^{-1}b_m$, then the topological group
$(G,\tau_1\lor\tau_2)$ is not complete and thus
is not determined by a $T$-sequence.
\endproclaim

\demo{Proof}
To prove the theorm it suffices to
verify that the diagonal of the square $G\ti G$ (which is identified with
$(G,\tau_1\lor\tau_2)$~) is not closed in the product
topology $\tau_1\ti\tau_2$.
Suppose the contrary: $G$ is closed in $G\ti G$. Then for the point
$(g,0)\in G\ti G$ beyond the diagonal we may find two neighborhoods
$U_i\in\tau_i$, $i=1,2$, of zero in $G$ such that $((U_1g)\ti
U_2)\cap G=\varnothing$. Since the sequences $(a_n)_{n\in\w}$ and
$(b_n)_{n\in\w}$ converge to zero in the topologies $\tau_1$,
$\tau_2$, we may find numbers $n,m\in\w$ such that $a_n\in U_1$,
$b_m\in U_2$, and $g=a_n^{-1}b_m$. Then $(b_m,b_m)=(a_ng,b_m)\in
(U_1g\ti U_2)\cap G$, a contradiction.
\qed
\enddemo

\remark{Remark 5}
Using Theorem 4 of \cite{Ba}, one may show that every non-closed
subgroup of a countable $k_\w$-group is not sequential. In particular,
the group $(G,\tau_1\lor\tau_2)$ constructed in
Theorem 5 is not sequential. Nonetheless, this group contains a non-trivial
convergent sequence (this can be easily shown using the sequentiality
of any countable $k_\w$-group, see Exercise 4.3.1 of \cite{PZ$_1$}).
\endremark
\vskip3pt

Finally, we prove the following theorem answering Question 2.5.5 of
\cite{PZ$_1$}.

\proclaim{Theorem 6} There exists a countable Abelian group $G$ admitting
a group topology $\tau$ determined by a $T$-sequence and a metrizable
group topology $\tau'$ such that the group $(G,\tau\lor\tau')$ is not
discrete but contains no non-trivial convergent sequence.
\endproclaim

\demo{Proof} Let $(G,\tau)=(\IZ^\w_f|(e_n))$ be the Graev
free topological Abelian group from Theorem 3. As we said the topology
$\tau$ is inductive with respect to the collection
$\{\IZ^\w_n\}_{n\in\w}$, where
$\IZ^\w_n=\{(x_i)_{i\in\w}:\sum_{i\in\w}|x_i|\le n\}$ for $n\in\w$.

Consider the metrizable topology $\tau'$ on $\IZ^\w_f$ generated by
the base $(U_n)_{n\in\w}$, where
$$
U_n=\{(x_i)_{i\in\w}\in\IZ^\w_f: x_i\in 2^n\cdot \IZ\text{ \ for all \ }
i\in\w\}, \;\; n\in\w.
$$
Let us show that the topology $\tau\lor\tau'$ is not discrete.
To see this, notice that for every $n\in\w$ and every open
neighborhood $U\in\tau$ of zero in $\IZ^\w_f$ the intersection $U_n\cap U$
is infinite (it contains the sequence $(2^ne_i)_{i\ge n_0}$ for some
$n_0$). Next, assume that $(b_n)_{n\in\w}$ is a sequence convergent to
zero in the topology $\tau\lor\tau'$.
Since $(b_n)_{n\in\w}$ is convergent in $(\IZ^\w_f,\tau)$
we get $\{b_n:n\in\w\}\subset \IZ^\w_{n_0}$
for some $n_0$. On the other hand, using the convergence of
$(b_n)_{n\in\w}$ to zero in $(\IZ^\w_f,\tau')$ we may find $m_0\in\IN$ such that
$\{b_n:n\ge m_0\}\subset U_{n_0}$. Since $U_{n_0}\cap Z^\w_{n_0}=\{0\}$,
we conclude that $b_n=0$ for all $n\ge m_0$, i.e., the sequence
$(b_n)_{n\in\w}$ is trivial.
\qed
\enddemo
\vskip5pt

{\bf Acknowledgement.} The author expresses his thanks to I.Protasov for
valuable and stimulating discussions on the subject of the paper.

\enddocument